\documentclass[11pt]{amsart}
\usepackage{verbatim, latexsym, amssymb, amsmath}
\usepackage{epsfig}
\def\R{\mathbb R}
\def\N{\mathbb N}
\def\Z{\mathbb Z}
\def\C{\mathbb C}

\def\H{\mathcal H}
\def\L{\mathcal L}
\def\x{\mathbf x}

\def\r{|x|}
\def\tr{|\partial_r^{\top}|}

\def\vol{\mathrm{vol}}
\def\re{\mathrm{Re}}
\def\hess{\mathrm{Hess}\,}

\def\d{\mathrm{div}}

\def\p{\mathrm{Proj}}

\newtheorem*{thma}{Theorem A}

\newtheorem{thm}{Theorem}[section]
\newtheorem{lemm}[thm]{Lemma}

\newtheorem{prop}[thm]{Proposition}
\theoremstyle{remark}
\newtheorem{rmk}[thm]{Remark}
\newtheorem{example}[thm]{Example}
\theoremstyle{definition}
\newtheorem{defi}[thm]{Definition}

\title{Translating solutions to Lagrangian mean curvature flow}

\author{{Andr\'e Neves} ${}^{\dagger}$}
\email{aneves@math.princeton.edu} 
\address{Fine Hall, Princeton University, Princeton, NJ 08544, USA}

\author{Gang Tian}
\email{tian@math.princeton.edu}
\address{Fine Hall, Princeton University,
Princeton, NJ 08544, USA}

\thanks{\quad\ ${}^{\dagger}$\ The author was partially supported by NSF grant DMS-06-04164.}

\begin{document}

\begin{abstract}

 We prove some non-existence theorems for translating solutions to Lagrangian mean curvature flow. More precisely, we show that  translating solutions with an $L^2$ bound on the mean curvature are planes and that  almost-calibrated translating solutions which are static are also planes.  Recent work of D. Joyce, Y.-I. Lee, and M.-P. Tsui, shows that these conditions are optimal.

\end{abstract}
\maketitle \markboth{Translating solutions to Lagrangian mean curvature flow} {Andr\'e Neves and Gang Tian}

\section{Introduction}\label{intro}

It was shown in \cite{neves} that finite-time singularities are, in some sense, unavoidable. More precisely, the first author gave examples of Lagrangians in $\C^2$ having the Lagrangian angle as small as we want and for which the Lagrangian mean curvature flow develops a finite-time singularity. Thus, if one aims to use Lagrangian mean curvature flow in order to understand the existence problem for Special Lagrangians (i.e. Lagrangians which are minimal surfaces) it is crucial to understand how finite-time singularities form.  The next example shows that this is a rather non-trivial problem.

\begin{figure}
\centering {\epsfig{file=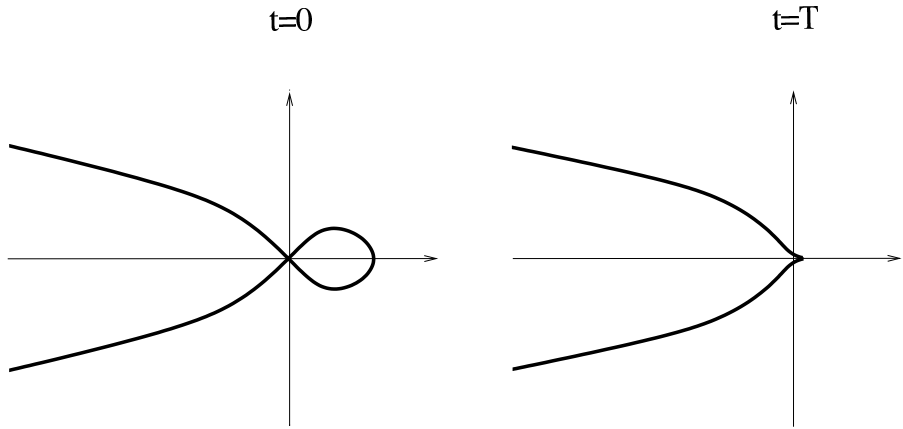, height=160pt}}\caption{Curve $\gamma_0$ and $\gamma_T$.} \label{fig2}
\end{figure}

\begin{example} \label{ex}Let $\gamma_0$ be the curve in $\C$ given in Figure \ref{fig2}.

 The curve can be made so that, under curve-shortening flow $(\gamma_t)_{t\geq 0}$, the small loop collapses at time $T$ and $\gamma_T$ becomes a curve with a cusp point. Moreover, $\gamma_0$ can be chosen so that the angle that the tangent vector makes with the $x$-axis has an oscillation not much bigger than $\pi$.
Let $L_t$ be the Lagrangian surface in $\C^2$ given by
$$L_t:=\gamma_t\times \R \subset \C\times\C.$$
Then, $L_0$ is a zero-Maslov class Lagrangian with oscillation of the Lagrangian angle as close to $\pi$ as we want and which develops a singularity at time $T$. The singular set is a line of cusp-points and hence has Hausdorff dimension one.
\end{example}
This example   shows that in order to develop a regularity theory for the flow, i.e., show that singularities for Lagrangian mean curvature flow are isolated, we need to require the oscillation of the Lagrangian angle to be strictly smaller than $\pi$ ({\em almost-calibrated}).

So far, the only known evidence that such regularity theory is possible was given in \cite{neves}. Their, assuming  the initial condition is a rational and almost-calibrated Lagrangian, the first author showed that if one rescales the flow around a  fixed point in space-time,  connected components of this rescaled flow converge to an area-minimizing union of planes.  The fact that the rescaled flow converges to a union of planes is an almost trivial consequence of Huisken's monotonicity formula and so the interesting part is that the configuration of planes needs to be area-minimizing. Without this property it would be hopeless to expect any regularity theory.

Nonetheless, we should point out that the property mentioned above is not sufficient  to develop a regularity theory.  One needs  to understand dilations of the flow where the point at which we center the dilation changes  with the scale (called {\em Type II dilations}).  From general theory, it follows that Type II dilations converge to an eternal solution with second fundamental form  uniformly bounded. If singularities are indeed isolated the expectation is that this eternal solution has vanishing mean curvature.
We now describe heuristically what could happen regarding Example \ref{ex}. If we rescale around the fixed point in space-time  at which the singularity is developing, the rescaled flow converges to a plane with multiplicity two. On the other hand, if we rescale the flow around the point of highest second fundamental form at some time $t_1$ close to $T$ and choose the scale so that the second fundamental form for the new solution becomes bounded by one, the rescaled flow  converges to the eternal solution given by
\begin{equation}\label{grim}
L_t:=\{(-\log \cos y_1+t, y_1, x_2,0)\,|\, -\pi/2< y_1 <\pi/2, \,x_2  \in \R\}.
\end{equation}
This solution is called the {\em grim-reaper} and is  an example of translating solutions to mean curvature flow.
\begin{defi}\label{trans}
A Lagrangian $L$ is a {\em translating solution} to Lagrangian mean curvature flow if we can find an ambient vector $e_1$ so that
\begin{equation*}
L_t:=L+te_1
\end{equation*}
is a solution to mean curvature flow.
\end{defi} 
\begin{rmk}
Without loss of generality, we assume that $e_1=(1,0,0,0)$. We can always achieve this by scaling $L$ and then  choosing a suitable coordinate system.
\end{rmk}

In a surprising new result, Dominic Joyce, Yng-ing Lee, and Mao-Pei Tsui \cite{JYT} found  translating solutions to Lagrangian mean curvature flow with oscillation of the Lagrangian angle arbitrarily small. They are described as follows. Let $w$ be a curve in $\C$ such that
$$w_t:=\sqrt{ 2t} w\quad\mbox{for}\quad t>0$$
is a solution to curve shortening flow in $\C$. This curve  can be chosen in a way that the angle $\theta$ that the tangent vector makes with the $x$-axis has arbitrarily small oscillation.
Set
\begin{equation}\label{jyt}
L:=\left\{\left(\frac{|w|^2(y)-x^2}{2}-i\theta(y), xw(y)\right),|\,x,y \in\R\right\}\subset \C\times \C.
\end{equation}
Using the fact that the curvature of $w$ satisfies
$$\vec k=w^{\bot},$$
it is a straightforward computation to check that $L$ is Lagrangian and that
$$L_t=L+t(1,0,0,0)$$
is a solution to Lagrangian mean curvature flow. Moreover, the Lagrangian angle of $L$ coincides with $\theta$ and hence its oscillation can be made arbitrarily small.

The main purpose of this paper is to give conditions that exclude 
 the existence of nontrivial translating solutions  to Lagrangian mean curvature flow. In order to do so, we need one more definition. 

Let $(L_t)_{-\infty<t<\infty}$ be a eternal solution to Lagrangian mean curvature flow in $\C^2$ which is almost-calibrated. Given a sequence $(\lambda_i)_{i\in\N}$ converging to zero, we can consider  the sequence of {\em blow-downs}
$$L^i_s:=\lambda_i L_{s/\lambda_i^2}\quad\mbox{where}\quad -\infty<s<\infty.$$ It follows from Theorem \ref{scale} that we can extract a subsequence  $L^i_s$ converging weakly to the same union of planes for every $s\leq 0$.

\begin{defi}
A eternal solution to  Lagrangian mean curvature flow is called {\em static} if we can find a  convergent sequence of blow-downs $L^i_s$  that converges weakly to the same union of planes for {\em every} $s$.
\end{defi}

In Theorem \ref{scale} we show that if $L_0$ is exact (see next section for definition), then $L^i_s$ converges to a self-expander for $s>0$ and if the self-expander has zero mean curvature, then the eternal solution is static.

\begin{thma}
Let $L$ be a translating solution for Lagrangian mean curvature flow in $\C^2$ for which we can find a constant $C_1$  such that
\begin{equation*}
		(H)\qquad\left \{\begin{aligned}
							&\mbox{The first Betti number of }L\mbox{ is finite,}\\
							&\sup_L |\theta|\leq C_1\\
							 &\H^2(L\cap B_R(0))\leq C_1 R^2\quad\mbox{for all }R>0,\\
							 &\sup_{L}|A|^2\leq C_1.
						\end{aligned}
					\right.
\end{equation*}
If 
$$\int_L |H|^2 d\mu \quad\mbox{is finite} $$
or if 
$$L \quad \mbox{is static}\quad\mbox{and}\quad\inf_L \cos \theta\geq\varepsilon_1$$
for some $\varepsilon_1>0$, then $L$ is a plane.
\end{thma}
 
The static condition is necessary for the following reason.  Consider the translating solution $(L_t)_{-\infty<t<\infty}$ described in \eqref{jyt} and denote by $\tilde w$ the curve in $\C$ given by the union of $w$ with $-w$. If $(\lambda_i)_{i\in\N}$ is a sequence converging to zero, the curves $\sqrt{\lambda_i}\tilde w$ converge to a union of two lines crossing at the origin which we denote by $\tilde w_0$. A simple computation shows that the sequence of blow-downs $L^i_s$ converges to the union of two planes given by
$$\R\times \tilde w_0\subseteq \C\times\C$$
for every $s\leq 0$ and  to a self-expander given by
$$\sqrt{2s}\left(\R\times\tilde w\right)=\R\times \sqrt{2s}\tilde w\subseteq\C\times\C$$
for every $s>0$.

The condition $$\inf_L \cos \theta\geq\varepsilon_1$$ is necessary because otherwise the grim-reaper described in \eqref{grim} would be a counterexample.

 The paper is organized as follows. In Section \ref{basic-defi} we introduce some basic notation and derive some simple identities for translating solutions and in Section \ref{auxi} we prove a compactness theorem for blow-down sequence of a eternal solution. 
 
 In Section \ref{mainsection} we assume that,  outside a compact set, $L$ can be decomposed into $N$ components $L_1,\ldots, L_N$, where each $L_j$ is the graph  of a multivalued function defined over a plane $P_j$ minus a disc, and show that in this case $L$ needs to be a plane  (Theorem \ref{main}).  The argument consists in  using  barriers to show that on each component $L_j$ the Lagrangian angle converges to a constant sufficiently fast and this will imply  that,  by choosing $R$ sufficiently large,  $L_j\setminus B_R(0)$ can be made as close to a plane as we want. In particular,
$$\lim_{R\to\infty}\oint_{\partial (L\cap B_R(0))}\langle\nu,e_1\rangle d\sigma=0,$$
	where $\nu$ denotes the exterior unit normal to $\partial (L\cap B_R(0))$ in $L$. The result will follow  because Proposition \ref{equacoes} (ii) implies that
	$$\oint_{\partial (L\cap B_R(0))}\langle\nu,e_1\rangle d\sigma=\int_{L\cap B_R(0)}|H|^2d\mu$$
	and so $L$ will have zero mean curvature.  In  that section we also prove a lemma (Lemma \ref{decomp2}) that gives us conditions under which a translating solution to Lagrangian mean curvature flow admits a graphical decomposition.

	In Section  \ref{final} we show that if $L$ has the $L^2$-norm of $|H|$ bounded  or if $L$ is static and almost-calibrated, then $L$ satisfies the conditions specified in Lemma \ref{decomp2} and hence admits a graphical decomposition.
	
\subsection{Open questions}

We propose two questions whose answer could provide some valuable insight on whether it is reasonable to expect any good behavior for singularities of Lagrangian mean curvature flow.

 The first question is whether the translating solutions found by  D. Joyce, Y.-I. Lee, and M.-P. Tsui, can arise as a blow-up of a finite time singularity fort Lagrangian mean curvature flow.

 If one aims to develop a regularity theory for the flow, it is absolutely necessary to answer this question.  This relation has been been observed before in codimension one mean curvature flow and Ricci flow. In the first case,   Brian White used his  work on mean convex solutions to mean curvature flow \cite{white1} to show that no grim-reaper appears as the limit of a sequence of rescaled flows \cite[Corollary 4]{white2}. In the second case, one of the first breakthroughs of Perelman  \cite[Section 4]{perelman} was to show that the cigar soliton does not arise as a finite-time singularity model. 
 
If one could show that, assuming the initial condition for the flow is an exact and almost-calibrated Lagrangian, no blow-down of a Type II rescale can give rise to a self expander, then the question above would be solved for a large class of initial conditions.

The second question addresses   the issue of uniqueness of translating solutions.  Suppose that $(L_t)_{t\in\R}$ and $(L'_t)_{t\in\R}$ are two translating solutions which are almost calibrated and exact. If after blow-down they produce the same self-expander, do they have to differ only by a rigid motion?
$$ $$

\noindent {\bf Acknowledgements:} The first author would like to express his gratitude to Yng-Ing Lee for some very stimulating  discussions.

\section{Definitions and basic identities} \label{basic-defi}
		
	Let $J$ and $\omega$ denote, respectively, the standard complex structure on $\C^2$ and the standard symplectic
	form on $\C^2$. We consider also the closed complex-valued $2$-form given by
	$$\Omega\equiv dz_1\wedge dz_2$$
	where $z_j=x_j+iy_j$ are complex coordinates of $\C^2$, and the Liouville form
	$$\lambda=\sum_{j=1}^2 x_jdy_j-y_jdx_j.$$
	We denote by $\H^1$ and $\H^2$ the one dimensional and two dimensional Hausdorff measures in $\C^2$ respectively.
	
	A smooth $2$-dimensional submanifold $L$ in $\C^2$ is said to be {\em Lagrangian} if $\omega_L=0$ and this
	implies that
	$$\Omega_L=e^{i\theta}\vol_L,$$
	where $\vol_L$ denotes the volume form of $L$ and $\theta$ is a multivalued function called the {\em
	Lagrangian angle}. When the Lagrangian angle is a single valued function the Lagrangian is called {\em
	zero-Maslov class} and if
	$$\cos \theta\geq \varepsilon_0$$ for some positive $\varepsilon_0$, then $L$ is said to be {\em almost-calibrated}.
	Furthermore, if $\theta\equiv\theta_0$, then $L$ is calibrated by
	$$\re\,\left( e^{-i\theta_0}\Omega\right)$$ and hence area-minimizing. In this case, $L$ is referred as being
	{\em Special Lagrangian}. The Lagrangian $L$ is said to be {\em exact} if the Liouville form is an exact form on $L$.
	
	Finally,  the
	relation between the Lagrangian angle and the mean curvature is given by 
	$$H=J\nabla \theta.$$

	Given a point $x_0$ in $\C^2$ and a time $T$, the backwards heat kernel is defined as
	$$\Phi_{x_0,T}(x,t)=\frac{\exp\left(-\frac{|x-x_0|^2}{4(T-t)}\right)}{4\pi(T-t)}.$$
	When $x_0$ is the origin and $T=0$, we denote it by $\Phi$. When it is clear from the context at which instant $t$ we are evaluating $\Phi_{x_0,T}(x,t)$, we denote it simply by $\Phi_{x_0,T}$. Moreover, $\x$ and $\x^{\bot}$ stand for the position vector associated with the point $x$ in $\C^2$ and its projection on the normal space of $T_x L$ respectively.
	
	Throughout this paper, $(L_t)_{-\infty< t <\infty}$ will be a translating solution   to Lagrangian mean curvature flow in $\C^2$, where $L$ satisfies hypothesis (H).  We denote by $|e_1^{\bot}|$ the projection of $e_1=(1,0,0,0)$ on the normal space of $L$. The intrinsic ball of radius $r$ around a  point $x$ in $L$ is defined by $\widehat B_r(x)$ and we fix $r_0<1$ to be such that $\widehat B_{r_0}(x)$ is simply connected for every $x$ in $L$. 
	
	\begin{prop}\label{equacoes}
		The following equations hold on $L$.
		\begin{itemize}
			\item[(i)]  There is a constant $\theta_0$ such that $$\theta(x)=-\langle J e_1,\x\rangle+\theta_0\qquad\mbox{and}\qquad |H|=|e_1^{\bot}|;$$
			\item[(ii)] $$\Delta x_1=|H|^2;$$
			\item[(iii)]$$\Delta \theta+\langle\nabla\theta,e_1\rangle=0.$$
		\end{itemize}
	\end{prop}
	\begin{proof}
		Because $(L_t)_{-\infty< t <\infty}$ is a translating solution to mean curvature flow we have that
		$$H=e_1^{\bot}$$\
		and thus
		$$\nabla \theta=-J e_1^{\bot}=-(Je_1)^{\top}.$$
		This implies the first property. The second one follows from
		$$\Delta x_1=\langle H, e_1\rangle=|H|^2$$
		and the third property is a consequence of
		$$\Delta \theta=-\d (Je_1)^{\top}=-\langle H, Je_1 \rangle=-\langle \nabla \theta, e_1\rangle.$$ 
	\end{proof}
	
	Given a sequence $(\lambda_i)_{i\in\N}$ converging to zero, we define the sequence of blow-downs
	$$L^i_s:=\lambda_i L_{s/\lambda_i^2}\quad\mbox{where}\quad -\infty<s<\infty.$$
	Finally, we we use the notation
	$$\{\theta=\alpha\}:=\{x\in L\,|\, \theta(x)=\alpha\}$$
	and
	$$\{|\theta-\alpha|\leq \delta\}:=\{x\in L\,|\,|\theta(x)-\alpha|\leq \delta\}.$$

\section{Blow-down Theorem}\label{auxi}

	Let $(\Sigma_t)_{-\infty<t<\infty}$ be an eternal solution to Lagrangian mean curvature flow for which we can find a constant $D$ such that, for every $t$,
	\begin{equation*}
		\qquad\left \{\begin{aligned}
							&  \H^2(\Sigma_t\cap B_R(0))\leq D R^2\quad\mbox{for all }R>0,\\
							&\cos(\theta_t)\geq D^{-1}.\\
						\end{aligned}
					\right.
	\end{equation*}	
	Given a sequence $(\lambda_i)_{i\in\N}$ converging to zero, consider the blow-downs
	$$\Sigma^i_{s}:=\lambda_i\Sigma_{s/\lambda_i^{2}}\quad\mbox{for}\quad{-\infty<s<\infty}.$$
	\begin{thm}\label{scale}
		There exist a finite set $$\{\bar\theta_1,\ldots,\bar\theta_N\}$$ and 
		Lagrangian planes
		$$P_1,\ldots,P_N$$
		such that, after passing to a subsequence, we have for every smooth function $\phi$ compactly supported, every
		$f$ in $C^2(\R)$, and every $s \leq 0$
		\begin{equation}\label{limite}
			\lim_{i \to \infty}\int_{\Sigma^i_s}f(\theta_{i,s})\phi d\mu=\sum_{j=1}^N m_j f(\bar\theta_j)\mu_j(\phi),
		\end{equation}
		where $\mu_j$ and $m_j$ denote the Radon measure of the support of $P_j$ and its multiplicity respectively. The set $$\{\bar\theta_1,\ldots,\bar\theta_N\}$$ does not depend on the sequence of rescales chosen.

		If $\Sigma_0$ is exact, there exists an integer rectifiable $2$-varifold  $\Sigma^{\infty}_1$ satisfying
		$$H=\frac{\x^{\bot}}{2}$$
		such that, after passing to a subsequence, $\Sigma^i_s$ converges as Radon measures to $\sqrt{s} \Sigma^{\infty}_1$ for every $s>0$.

		If $\Sigma^{\infty}_1$ is stationary, then identity \eqref{limite} holds for  {\em every} s and hence $(\Sigma_t)$ is static.
	\end{thm}
	
	\begin{rmk}
	\begin{itemize}
		\item [(i)]Whenever identity \eqref{limite} holds we say that $\Sigma^i_s$ converges weakly to $\Sigma^{\infty}_0$.
		\item[(ii)] In Section \ref{intro} we saw that the last property of the theorem does not hold when $L$ is the grim-reaper in $\C^2$. Hence, we see that the almost-calibrated condition is necessary in order for the convergence to hold when $s=0$.
	\end{itemize}
	\end{rmk}
	
	\begin{proof}
		
	The first property was essentially  proven in \cite[Theorem A]{neves}. The ideas apply with no modification.
	
	From the compactness for integral Brakke motions \cite[Section 7.1]{ilmanen1} we know that, after passing to a subsequence, $(\Sigma^i_s)_{-\infty<s<\infty}$ converges to an integral Brakke motion $(\Sigma^{\infty}_s)_{-\infty<s<\infty}$. Moreover, Federer and Fleming compactness for integral currents implies that $\Sigma^i_0$ converges to a current $\Sigma.$ The fact that $\Sigma$ is almost-calibrated implies that the support of $\Sigma$ equals the support of $\Sigma^{\infty}_0$ because, for every $\phi\geq 0$,
			$$
			\int_{\Sigma^{\infty}_0}\phi\, d\mu \leq D\lim_{i\to\infty}\int_{\Sigma^i_0}\phi\cos\theta\, d\mu=D\int_{\Sigma} \phi \mbox{Re\,}\Omega.
			$$
			From the monotonicity formula \cite{huisken} we have that, for every $x_0$ in $\C^2$, every positive $T$, and  every $s<0$
			$$\int_{\Sigma^{i}_0}\Phi_{x_0, T}(\cdot, 0)\, d\mu\leq \int_{\Sigma^{i}_s}\Phi_{x_0, T}(\cdot, s)\, d\mu$$
			and thus
			\begin{equation*}\label{monotone}
				\int_{\Sigma^{\infty}_0}\Phi_{x_0, T}(\cdot, 0)\, d\mu\leq  \sum_{j=1}^N m_j \int_{P_j}\Phi_{x_0, T}(\cdot, s)\,d\mu.
			\end{equation*}
			Denote the density function of $\Sigma^{\infty}_0$ and $\Sigma^{\infty}_{-1}$ by $\Theta_0(x)$ and $\Theta(x)$ respectively. Make $T$ go to zero so that the left-hand side converges to $\Theta_0(x_0)$ for almost all $x_0$. After that, make $s$ converge to zero to obtain  $\Theta(x_0)\geq \Theta_0(x_0)$ and so the support of $\Sigma$ is contained on the support of $\Sigma^{\infty}{-1}$. Because $\partial \Sigma=0$, the Constancy Theorem \cite[Theorem 26.27]{Leon} implies that the support of $\Sigma$ coincides with $\Sigma^{\infty}{-1}$. We now argue that the density functions also coincide.

			 We know that
			$$\frac{d \cos \theta_t }{dt}=\Delta \cos \theta_t +\cos\theta_t|H|^2$$
			and thus, for every $T>0$, $x_0$ in $\C^2$, and $s<0$, we have from the monotonicity formula that
			\begin{multline*}
				\frac{d}{ds}\int_{\Sigma^i_s}  \cos \theta_{i,s} \Phi_{x_0,T}\,d\mu \\
				=\int_{\Sigma^i_s} \cos \theta_{i,s}\left(|H|^2-\left|H+\frac{(\x-\x_0)^{\bot}}{2(T-s)}\right|^2\right)\Phi_{x_0,T}\,d\mu\\
				\geq -\int_{\Sigma^i_s} \cos \theta_{i,s}\left(\delta |H|^2+C\left|\frac{(\x-\x_0)^{\bot}}{2(T-s)}\right|^2\right)\Phi_{x_0,T}\,d\mu
			\end{multline*}
			where $C=C(\delta)$. The evolution equation for $\theta_{i,s}$ implies that for all $s<0$
			\begin{multline*}
			\lim_{i\to\infty}\int_s^0 \int_{\Sigma^i_t} 2 |H|^2 \Phi_{x_0,T}d\mu\,dt\leq \int_{\Sigma^i_s}  \theta^s_{i,s} \Phi_{x_0,T}d\mu\\
			 =  \sum_{j=1}^N m_j \bar\theta_j^2 \int_{P_j} \Phi_{x_0,T}(\cdot,s)\,d\mu\leq \sum_{j=1}^N m_j \bar\theta_j^2 :=B
			\end{multline*}
			and it is not hard to see that
			$$\lim_{i\to\infty}\int_s^0 \int_{\Sigma^i_t} \left|\frac{(\x-\x_0)^{\bot}}{2(T-t)}\right|^2 \Phi_{x_0,T}d\mu\,dt=0.$$
			Thus
			\begin{multline*}
				\int_{\Sigma} \Phi_{x_0,T}(\cdot,0)\,\mbox{Re}\,\Omega  =\lim_{i\to\infty}\int_{\Sigma^i_0} \cos \theta_{i,0} \Phi_{x_0,T}\,d\mu\\
				\geq -\delta B+\lim_{i\to\infty} \int_{\Sigma^i_s} \cos(\theta_{i,s})\Phi_{x_0,T}\,d\mu\\
				=-\delta B+\sum_{j=1}^N m_j \cos\bar\theta_j\int_{P_j} \Phi_{x_0,T}(\cdot,s)\,d\mu.
			\end{multline*}
			Making $\delta,$ $T$, and $s$ converge to zero, and using the almost-calibrated condition,  we obtain that for almost all $x_0$
			$$ \Theta_0(x_0)\geq \Theta(x_0)$$
			and so $\Sigma$ coincides with $\Sigma^{\infty}_{-1}$.
			As a result, we have that for every $T$ positive
			$$\lim_{i\to\infty}\int_{\Sigma^i_{-1}}\Phi_{0,T}(\cdot,-1)d\mu=\lim_{i\to\infty}\int_{\Sigma^i_{0}}\Phi_{0,T}(\cdot,0)d\mu$$
			and so the monotonicity formula implies that
			$$\lim_{i\to\infty}\int_{-1}^0 \int_{\Sigma_i^s} (|H|^2+|\x^{\bot}|^2)\exp(-|x|^2)d\mu ds=0.$$
			We can then argue as in \cite[Theorem A]{neves} and conclude that identity \eqref{limite} holds for $s=0$.

			Recall that $\Sigma^i_0$ converges as Radon measures and as currents to
	a union of planes with possible multiplicities.  Because the Grassmanian of $2$-planes in $\R^4$ can be parametrized by  the self-dual and anti self-dual two forms of $\R^4$, we have that
	$$\lim_{i\to\infty} \int_{\Sigma^i_0\cap B_R(0)} |\x^{\bot}|^2d\mu=0$$
	for every positive $R$. Moreover, the fact that $\Sigma_0$ is exact implies the existence of $\beta^i_s$ defined on $\Sigma^i_s$ for which
	$$d\beta^i_s=\lambda\quad\mbox{and}\quad |\nabla \beta^i_s|=|\x^{\bot}|.$$
	Hence, because $\Sigma^i_0$ is connected, we can apply Proposition \ref{poincare} and conclude the existence of $\bar \beta$ such that, after passing to a subsequence,
	$$\lim_{i\to\infty}\int_{\Sigma^i_0\cap B_R(0)}(\beta^i_0-\bar\beta)^2d\mu=0.$$
	According to \cite[Section 6]{neves}, we have that
	$$\frac{d}{ds}(\beta^i_s+2s\theta^i_s)^2=\Delta(\beta^i_s+2s\theta^i_s)^2-2|\x^{\bot}-2sH|^2$$
	and so we can apply Huisken's monotonicity formula to conclude that, for every positive $T$,
	\begin{equation}\label{expand}
	\lim_{i\to\infty}\int_0^T\int_{\Sigma^i_s} 2|\x^{\bot}-2sH|^2 \Phi_{0,T}d\mu ds\leq\lim_{i\to\infty} \int_{\Sigma^i_0} (\beta^i_0-\bar\beta)^2\Phi_{0,T} d\mu=0.
	\end{equation}
	We want to show that, for every compactly supported function $\phi$ in $\C^2$,
	$$\frac{1}{s}\int_{\Sigma^{\infty}_s} \phi\left({x}/{\sqrt s}\right)d\mu$$
	is constant as a function of $s$ for every $s>0$.
	
	A standard computation shows that
	\begin{multline*}
	\frac{d}{ds} \left (\frac{1}{s}\int_{\Sigma^{\infty}_s} \phi\left({x}/{\sqrt s}\right)d\mu \right) =-\frac{1}{s^2}\int_{\Sigma^{\infty}_s} \phi d\mu-\frac{1}{2s^{5/2}}\int_{\Sigma^{\infty}_s} \langle D\phi, \x \rangle d\mu\\+
	\frac{1}{s^{3/2}}\int_{\Sigma^{\infty}_s} \langle D\phi, H \rangle d\mu-\frac{1}{s}\int_{\Sigma^{\infty}_s} \phi|H|^2 d\mu.
 	\end{multline*}
	Due to
	$$\Delta |x|^2=4+2\langle\x,H\rangle$$
	we obtain that
	$$2\int_{\Sigma^i_s}\phi\left({x}/{\sqrt s}\right)d\mu=-\int_{\Sigma^i_s}\langle \x,H\rangle\phi d\mu-\frac{1}{\sqrt s}\int_{\Sigma^i_s}\langle \x^{\top}, D\phi\rangle d\mu.$$
	Hence,
	\begin{multline*}
	\frac{d}{ds} \left (\frac{1}{s}\int_{\Sigma^{\infty}_s} \phi\left({x}/{\sqrt s}\right)d\mu \right) =\frac{1}{s^{3/2}}\int_{\Sigma^{\infty}_s}\left \langle D\phi, H-\frac{x^{\bot}}{2s} \right\rangle d\mu \\+\frac{1}{s}\int_{\Sigma^{\infty}_s}\phi \left \langle H, \frac{x^{\bot}}{2s}-H \right\rangle d\mu.
	\end{multline*}
	For every $0<a<b$ and $R>0$, we have that
	$$\int_a^b \int_{\Sigma^i_s}|D\phi|^2d\mu ds\quad\mbox{and}\quad\int_a^b \int_{\Sigma^i_s\cap B_R(0)}|H|^2d\mu ds$$
	are uniformly bounded. Therefore, \eqref{expand} implies that
	$$\frac{1}{s}\int_{\Sigma^{\infty}_s} \phi\left({x}/{\sqrt s}\right)d\mu$$
	is indeed independent of $s$ for all $s>0$ and this is equivalent to $\Sigma^{\infty}_s=\sqrt s \Sigma^{\infty}_1$ for all $s>0$.

			We now show the last property. Because $\Sigma^{\infty}_1$ has $\x^{\bot}=0$, we obtain that for every $0<a<b$ and every $R>0$
			$$\lim_{i\to\infty}\int_{a}^b\int_{\Sigma^i_s\cap B_R(0)}|\x^{\bot}|^2d\mu ds.$$
			Hence, identity \eqref{expand} implies that
			$$\lim_{i\to\infty}\int_{a}^b \int_{\Sigma_s^i} (|H|^2+|\x^{\bot}|^2)\exp(-|x|^2)d\mu ds=0$$
			and so, assuming without loss of generality that
			$$\int_{\Sigma_1^i} (|H|^2+|\x^{\bot}|^2)\exp(-|x|^2)d\mu ds=0,$$
			\cite[Proposition 5.1]{neves} implies that $\Sigma^{\infty}_1$ is a union of Lagrangian planes with multiplicities. In order to prove the result, it suffices to show that
			$\Sigma^{\infty}_1$ equals $\Sigma^{\infty}_{0}$.
			
			 Huisken's monotonicity formula implies that for every $T>0$ and $0<s<T$
			$$\int_{\Sigma^{\infty}_0}\Phi_{x_0,T} d\mu\leq \int_{\Sigma^{\infty}_s}\Phi_{x_0,T}(\cdot,s)d\mu=\int_{\Sigma^{\infty}_1}\Phi_{x_0,T}(\cdot,s)d\mu.$$
			Thus making $s$ converge to $T$ and then $T$ converge to zero, we obtain that
			$$\Theta_{0}(x_0)\geq \Theta_1(x_0),$$
			where $\Theta_1$ denotes the density function of $\Sigma^{\infty}_1$. Hence the support of $\Sigma^{\infty}_1$ is contained in the support of $\Sigma^{\infty}_{0}$.
		
			Choosing $x_0$ to be in exactly one $P_{j_0}$, we can apply Huisken's monotonicity formula to $\cos \theta^i_s$ and conclude that for any  $0<s<T$
			\begin{multline*}
				\int_{\Sigma^{\infty}_1}\cos(\bar\theta_j) \Phi_{x_0,T}(\cdot,s)\,d\mu=\lim_{i\to\infty}\int_{\Sigma^i_s} \cos(\theta_{i,s}) \Phi_{x_0,T}(\cdot,s)\,d\mu\\
			\geq \lim_{i\to\infty }\int_{\Sigma^i_0} \cos(\theta_{i,s}) \Phi_{x_0,T}(\cdot,0)\,d\mu =
			\sum_{j=1}^N m_j \cos(\bar\theta_j)\int_{P_j} \Phi_{x_0,T}(\cdot,0)\,d\mu\\
			\geq \cos(\bar\theta_{j_0})m_{j_0}.
			\end{multline*}
			Making $s$ converge to $T$, the almost-calibrated condition  implies that  $$\Theta_{0}(x_0)\leq \Theta_1(x_0).$$			
	\end{proof}

\section{Graphical implies flatness}\label{mainsection}

We start by defining what it means for a Lagrangian $L$ to have a graphical decomposition. Recall that $r_0$ was chosen (see Section \ref{basic-defi}) so that  $\widehat B_{r_0}(x)$ is simply-connected for all $x$ in $L$.

\begin{defi}\label{GD}A Lagrangian $L$ is said to admit a {\em graphical decomposition} if,  outside a compact set, $L$ can be decomposed into $N$ connected components $L_j$, $j=1,\ldots,N,$ having the following property.
				
	For each $j=1,\ldots,N$ there are constants $\bar \theta_j$,  $R_j$, $S_j$, $n_j$ positive integer, and  a Lagrangian plane $P_j$
			$$P_j:=\{(u,0,v\cos \bar \theta_j ,v\sin \bar \theta_j)\,|\, (u,v)\in \R^2\}$$
	such that
	\begin{itemize}
		\item[(i)] For every $x$ in $L_j$, $\widehat B_{s_0}(x)$ can be written as the graph of a function defined over $P_j$ with its derivatives  bounded by $S_j$;		
		\item[(ii)] The projection of $L_j$ on $P_j$ $$\p_{P_j}\colon L_j \longrightarrow P_j\setminus B_{R_j}(0)$$ is a $n_j$-fold covering map;			
		\item[(iii)] 
				$$\lim_{R\to\infty} \sup_{L_j\setminus B_R(0)}|\theta-\bar\theta_j|=0.$$
	\end{itemize}
	In particular,
	$$\lim_{R\to\infty}\sup_{L\setminus B_R(0)}|H|=0. $$
\end{defi}
\begin{rmk}
		In case $L$ is  not embedded, property (ii) in Definition \ref{GD} should be interpreted as follows. If $F$ denotes the immersion of the surface $L$ in $\C^2$, then $\p_{P_j}\circ F$ is  a $n_j$-fold covering map.
\end{rmk}

In this section we show that any translating solution  $L$ to Lagrangian mean curvature flow which admits a graphical decomposition is a plane.

\begin{thm}\label{main}If $L$ is a translating solution to Lagrangian mean curvature flow with uniformly bounded second fundamental form and admitting a graphical decomposition, then $L$ is a plane.
\end{thm}
\begin{proof}
It suffices to show that the mean curvature of $L$ is zero because, in that case, we have from  Proposition \ref{equacoes} that $L$  is a Special Lagrangian contained in some hyperplane and therefore a plane. The idea consists in finding a sequence of compact sets $K_n$ exhausting $L$ for which
	$$\lim_{n\to\infty}\oint_{\partial K_n}\langle\nu,e_1\rangle d\sigma=0,$$
	where $\nu$ denotes the exterior unit normal to $\partial K_n$ in $L$. The theorem follows because, due to Proposition \ref{equacoes} (ii),
	$$\oint_{\partial K_n}\langle\nu,e_1\rangle d\sigma=\int_{K_n}|H|^2d\mu.$$ 

The following barrier will be needed to prove the main theorem.

\begin{lemm}\label{barrier}
		For every $\alpha<1/3$, there is a constant $R_0=R_0(\alpha)$ such that, for every constants $\delta$ and $B$, the function
		$$V_{\delta, B}:=B\r^{-\alpha}\exp(-\r/2)+\delta\exp(x_1/2)$$
		satisfies
		$$\Delta V_{\delta,B}\leq V_{\delta,B}\frac{|H|^2+1}{4}\quad\mbox{for all}\quad \r\geq R_0.$$
		\end{lemm}
	\begin{proof}
		Denote by $\partial_r$, $\partial_r^{\top}$, and $\partial_r^{\bot}$ the radial vector,  its tangential projection on $L$, and its projection on the normal bundle of $L$ respectively.  Set $$f(\r):=\r^{-\alpha}\exp(-\r/2).$$
		Away from the origin, we have
		$$\Delta \r=\langle \partial_r,H\rangle+\frac{2}{\r}-\frac{\tr^2}{\r}$$
		and so, because
		$$f'=-\alpha\frac{f}{\r}-\frac{f}{2}\quad\mbox{and}\quad f''=\alpha(\alpha+1)\frac{f}{\r^2}+\alpha\frac{f}{\r}+\frac{f}{4},$$
		we obtain that
		\begin{multline*}
		\Delta f=f\left(\frac{\tr^2}{4}-\frac{\langle \partial_r,H\rangle}{2}\right)+\frac{f}{|x|}\left(\alpha\tr^2-1-\alpha\langle \partial_r,H\rangle+\frac{\tr^2}{2}\right)\\
		+\frac{f}{|x|^2}\left(\alpha(\alpha+2)\tr^2-2\alpha\right).
		\end{multline*}
		Hence
		\begin{align*}
			\Delta f &\leq f\frac{|H|^2+1}{4}+ \frac{f}{|x|}\left(\frac{\alpha|H|^2}{2}+\alpha-\frac{1}{2}\right)+\frac{f}{|x|^2}{\alpha^2}\\
			&\leq f\frac{|H|^2+1}{4}+\frac{f}{|x|}\frac{3\alpha-1}{2}+\frac{f}{|x|^2}{\alpha^2}.
		\end{align*}
		Therefore, we can choose $R_0=R_0(\alpha)$ so that for all $\r\geq R_0$ we have
		$$\Delta f \leq f\frac{|H|^2+1}{4}.$$
		Proposition \ref{equacoes} (ii) implies that
		$$\Delta \exp(x_1/2) = \exp(x_1/2)\frac{|H|^2+1}{4}$$
		and so the proposition follows.
	\end{proof}

	This proposition implies the following decay for $|\theta-\bar\theta_j|$ on each component $L_j$, $j=1,\ldots,N$ given by the graphical decomposition.
	
	\begin{lemm}\label{decay} For each $\alpha<1/3$, there are constants $B$ and  $R_0$ such that, on each component $L_j$,  $j=1,\ldots,N$, we have
		$$|\theta(x)-\bar \theta_j|+|\nabla\theta|(x)\leq B\r^{-\alpha}\exp(-\r/2-x_1/2)\quad\mbox{for all}\quad \r\geq R_0.$$
	\end{lemm}
	\begin{proof}
		For each $j=1,\ldots,N,$ set
		$$u_j:=(\theta-\bar \theta_j)\exp(x_1/2).$$
		From Proposition \ref{equacoes} we know that
		$$\Delta u_j = u_j\frac{|H|^2+1}{4}$$
		and thus, we obtain from Lemma \ref{barrier}  that
		$$\Delta (V_{\delta,B}-u_j) \leq (V_{\delta,B}-u_j) \frac{|H|^2+1}{4}\quad\mbox{for all}\quad \r\geq R_0,$$
		where $R_0$ is chosen large enough so that  $\partial L_j\subseteq B_{R_0}(0)$ and 		
		the constant $B$ is chosen so that, for every $j=1,\ldots,N$, we have
		$$x\in L_j\cap\partial B_{R_0}(0) \implies B\r^{-\alpha}\exp(-\r/2)>|u_j(x)|.$$
		 Note that, for all $R$ sufficiently large, we have from Definition \ref{GD} (iii)
		 $$\sup_{L_j\cap\partial B_R(0)}|u_j\exp(-x_1/2)|\leq \delta$$	
		 and thus 
		 $$\sup_{L_j\cap\partial B_R(0)}(V_{\delta,B}-u_j)>0 \quad \mbox{for every}\quad j=1,\ldots,N.$$
		 Applying  the maximum principle  to $L_j\cap (B_R(0)\setminus B_{R_0}(0))$ for all $R$ sufficiently large, we have that
		 $$u_j(x)\leq B\r^{-\alpha}\exp(-\r/2)+\delta\exp(x_1/2) \quad\mbox{for all}\quad \r\geq R_0.$$ 
		 As a result, after making $\delta$ tend to zero, we obtain
		 $$\theta(x)-\bar \theta_j\leq B\r^{-\alpha}\exp(-\r/2-x_1/2) \quad\mbox{for all}\quad \r\geq R_0.$$
		 The correspondent estimate for $|\theta(x)-\bar \theta_j|$ follows in the same way by considering the function $V_{\delta,B}+u_j$. The estimate for $|\nabla \theta|$ for  is a consequence of Proposition \ref{equacoes} (i) and interior Schauder estimates.
	\end{proof}

	The graphical decomposition   implies the existence of $r_1$ so that, for every $p$ in $L_j$, $j=1,\ldots,N$ the projection of $\widehat B_{r_0}(p)$ on $P_j$ contains $B_{r_1}(\bar p)\cap P_j$, where $\bar p$ stands for the projection of $p$ on $P_j$. Thus, after an appropriate change of coordinates, a neighborhood of $p$ can be described as
	 $$(u,v,\partial_u f,\partial_v f)\quad\mbox{with}\quad (u,v)\in B_{r_1}(\bar p)$$
	 for some function $f$ with $|\hess f|$ uniformly bounded, where the coordinate $x_1$ equals $u$ and the coordinate $y_1$ equals   $\partial_u f$. Furthermore,  a direct computation shows that we can find a constant $D$ depending only on the constants $S_j$ (see Definition \ref{GD} (i)) for which one of the eigenvalues of $\hess f$ satisfies
	 $$|\lambda|\leq D |e_1^{\bot}| \quad\mbox{on}\quad B_{r_1}(\bar p).$$
	 Thus, we obtain from Lemma \ref{decay} that, provided $|p|$ is large enough, one of the eigenvalues $\lambda$ has the decay
	 \begin{equation}\label{eq}
	 	|\lambda| \leq  B\r^{-\alpha}\exp(-\r/2-u/2)
	 \end{equation}
	 for some constant $B$.
	 
	 An explicit computation shows that  the function $f$ satisfies the equation
	$$\arctan \lambda_1+\arctan \lambda_2=\theta(x)-\bar\theta_j,$$
	where $\lambda_1$ and $\lambda_2$ are the eigenvalues of $\hess f$.  As a result, Lemma \ref{decay} and estimate \eqref{eq} imply that, provided $|p|$ is large enough,
	\begin{equation}\label{hess}
		|\hess f| \leq  B\r^{-\alpha}\exp(-\r/2-u/2),
	\end{equation}
	where $\alpha<1/3$ and $B$ depends on $\alpha$.

	On each of the connected components $L_j$, denote by $\gamma_{j,r}$ the lift to $L_j$ of the path on $P_j$ given by
	$$c_j(t):=(r\cos t,r\sin t)\quad\quad 0\leq t\leq 2n_j\pi,$$
	where the variable $r$ will be made as large as we want.
	
	There is $t_0=t_0(r,r_1)$ such that, for every $ t_1\leq 2n_j\pi,$ we can find a function $f$ for which
	$$\gamma_{j,r}(t)=(r\cos t,r\sin t,\partial_u f, \partial_v f)\quad\mbox{for all } t_1-t_0\leq t \leq t_1+t_0,$$
	where, due to \eqref{hess}, the function $f$ restricted to this portion of $\gamma_{j,r}$ satisfies
	$$|\hess f|(t) \leq  B r^{-\alpha}\exp(-r/2-r\cos t/2).$$
	Using an obvious abuse of notation, the tangent vector $\gamma_{j,r}'(t)$ is given by
	$$\gamma_{j,r}'(t)=r(\partial_t, (\hess f)(\partial_t)),\quad\mbox{where}\quad \partial_t:=(-\sin t,\cos t)$$
	and, denoting by $\bar\nu$ the vector in $\R^2$ for which
	$$\nu=(\bar\nu, (\hess f)(\bar\nu)),$$
	then
	$$\langle \bar\nu,\partial_t+(\hess f)^{\ast}( \hess f) (\partial_t)\rangle=0\quad\mbox{and}\quad 1=|\bar\nu|^2+|(\hess f) (\bar\nu)|^2, $$
	where $(\hess f)^{\ast}$ denotes the transpose of $\hess f.$
	Thus, provided we choose $r$ sufficiently large, we can find a constant $C$ such that
	$$|\bar\nu-\partial_r|\leq  C r^{-2\alpha}\exp(-r-r\cos t)$$
	and 
	$$|\gamma_{j,r}'(t)|\leq r(1+ Cr^{-2\alpha}\exp(-r-r\cos t)).$$
	For this reason,
	$$\int_{t_1-t_0}^{t_1+t_0}\langle\nu,e_1 \rangle |\gamma_{j,r}'(t)|\,dt=\int_{t_1-t_0}^{t_1+t_0}r\cos t\,dt+Q,$$
	where we can find a constant $C$ such that $$|Q|\leq \int_{t_1-t_0}^{t_1+t_0} C r^{1-2\alpha}\exp(-r-r\cos t)\,dt.$$
	Therefore,  we obtain from unique continuation that
	\begin{align*}
		\left|\oint_{\gamma_{j,r}}\langle\nu,e_1\rangle d\sigma\right| & \leq\int_{0}^{2n_j\pi} C r^{1-2\alpha}\exp(-r-r\cos t)\,dt\\
		& = n_jC r^{1-2\alpha}\int_{0}^{2\pi}\exp(-r-r\cos t)\,dt.
	\end{align*}
	 Choose $\delta$ so that
	$$|\cos(y)+1|\leq (y-\pi)^2\quad\mbox{for all }|y-\pi|\leq \delta.$$
	Then, we can find a positive constant $D=D(\delta)$ so that, for all $r$ sufficiently large,
	\begin{align*}
		\int_{0}^{2\pi}\exp(-r-r\cos t)\,dt\leq & \int_{\pi-\delta}^{\pi+\delta}\exp(-r(t-\pi)^2)\,dt+2\pi \exp(-Dr)\\
		\leq & r^{-1/2}\int_{-\sqrt r\delta}^{\sqrt r \delta}\exp(-s^2)\,ds+2\pi \exp(-Dr)\\
		\leq &  Dr^{-1/2}.
	\end{align*}
	Hence, provided we choose $\alpha>1/4$ in Lemma \ref{decay}, we have that 
	$$\lim_{r\to\infty}\left|\oint_{\gamma_{j,r}}\langle\nu,e_1\rangle d\sigma\right|= \lim_{r\to\infty} r^{1/2-2\alpha}=0.$$
	Property (i) and (ii) of Definition \ref{GD} imply that we can find a sequence of compact sets $K_n$ exhausting $L$ and such that
	$$\partial K_n=\gamma_{1,r_n}\cup \cdots \cup \gamma_{N,r_n},$$
	where $(r_n)_{n\in\N}$ is a sequence converging to infinity. This finishes the proof.

\end{proof}

The next lemma gives conditions under which a translating solution $L$  satisfying (H) admits a graphical decomposition.

	\begin{lemm}\label{decomp2}
				Assume that for almost all $-C_1<a<C_1$ we have
				$$\{\theta=a\}\subset B_{R(a)}	$$
				for some positive $R(a)>0$. 				
				Then $L$ admits a graphical decomposition.			
	\end{lemm}
	\begin{proof}
		For every $\alpha$, there is only one Lagrangian plane $P_\alpha$ with Lagrangian angle $\alpha$ and $|e_1^{\bot}|=0$ which is given by
		\begin{equation}\label{plane}
			P_{\alpha}=\{(u,0,v\cos \alpha ,v\sin \alpha)\,|\, (u,v)\in \R^2\}.
		\end{equation}
		  Set $\omega_{\alpha}$ to be the volume form of that plane extended by parallel translation to $\C^2$ and denote its Hodge-dual on $L$ by $\ast \omega_{\alpha}$. 
		  
		  We claim that for every  $\varepsilon$, there is $\delta_1$ such that, for every $x$ in $L$ and every $\alpha$, we have
		  \begin{equation}\label{cont}
		  \sup_{\widehat B_{r_0}(x)}|\theta-\alpha|\leq \delta_1\quad\implies\inf_{\widehat B_{r_0}(x)}\ast\omega_{\alpha}\geq \varepsilon.
		  \end{equation}		  	
		  A simple continuity argument shows the existence of $\delta_2$ such that  if $Q$ is a Lagrangian plane with Lagrangian angle $\theta(Q)$, then for every $\alpha$
		$$|e_1^{\bot}|\leq \delta_2\quad\mbox{and}\quad |\theta(Q)-\alpha|\leq\delta_2\quad\implies\quad \omega_{\alpha}(Q)\geq \varepsilon,$$
		where $\omega_{\alpha}(Q)$ denotes the evaluation of $\omega_{\alpha}$ on $Q$.
		Moreover, as we argue next, we can find $\delta_1\leq \delta_2$ such that for all $x$ in $L$
		$$\sup_{\widehat B_{r_0}(x)}|\theta-\theta(x)|\leq \delta_1/2\quad\implies\quad \sup_{\widehat B_{r_0}(x)}|H|\leq \delta_2.$$
		If not, we could find a sequence of translating solutions $(L_j)_{j\in\N}$ converging smoothly on compact sets to another translating solution $L_{\infty}$ with Lagrangian angle constant on $\widehat B_{r_0}(0)$ and $|H|$ not identically zero on $\widehat B_{r_0}(0)$. This proves the desired claim.		
		  
		Let   $(x_i)_{i\in\N}$ be a sequence in $L$ with $|x_i|$ going to infinity and $\theta(x_i)$ converging to some $\alpha$. Set $\varepsilon=1/2$ in identity \eqref{cont} and choose $\delta<\delta_1/2$ so that
		$$\{\theta=\alpha\pm\delta\}\subset B_{R(\alpha\pm\delta)}.$$
		Set
		 $$R:=\max \{R(\alpha+\delta),R(\alpha-\delta)\}+2\quad\mbox{and}\quad\Sigma:=\{|\theta-\alpha|\leq \delta\}.$$
		
		We have
		$$\partial \Sigma \subseteq B_{R-2}(0)$$
		and if $x$ belongs to $\Sigma\setminus B_R(0)$, identity \eqref{cont} implies that
		$$\inf_{\widehat B_{r_0}(x)}\ast\omega_{\alpha}\geq 1/2.$$
		Moreover, from \cite[page 476]{neves}, there exists a constant $D$ such that
		\begin{equation}\label{density}
			\inf_{x\in L,\, r\leq r_0} r^{-2}\H^2(\widehat B_{r}(x))\geq D^{-1}
		\end{equation}
		and so we can apply Lemma \ref{basic2} in order to conclude that, outside a compact set $K$, 
		$$\Sigma\setminus K=L_1\cup\cdots\cup L_{N_1}$$
		where,  for each $j=1,\ldots, N_1,$ 
		$$\p_{P_{\alpha}}\colon \Sigma_j \longrightarrow P_{\alpha}\setminus B_{R}(0)$$ is a $n_j$-fold covering map and, for every $x$ in $\Sigma\setminus B_R(0)$, $\widehat B_{r_0}(x)$ can be written as the graph of a function defined over $P_{\alpha}$ with its derivatives uniformly bounded.
		
		Take a unbounded sequence $(y_i)$ in $L_j$ such that $\theta(y_i)$ converges to some $\bar \theta_j$  and suppose there is another unbounded sequence $p_i$ in $L_j$ such that $\theta(p_i)$ converges to some $\beta$ distinct from $\bar \theta_j$. Using the map $\p_{P_{\alpha}}$ and the local graphical property we can find, for every $a$  between $\beta$ and $\bar\theta_j$, an unbounded sequence $w_i$ in $L_j$ such that $\theta(w_i)=a$. This contradicts our hypothesis  and proves the third property of Definition \ref{GD}. The first and second property of Definition \ref{GD} follow because, by choosing $R_j>R$,  we can replace $\alpha$ by $\bar \theta_j$ on each $L_j$.
				
		We can repeat the whole process but this time applied to $L\setminus \Sigma$. We only need to this finitely many times because 
				$$\lim_{R\to\infty} \H^2(L_j\cap B_R(0))R^{-2}\geq C$$
		for some universal constant $C$. 
	\end{proof}

\section{Proof of main theorem}\label{final}

We now prove Theorem A. We start by showing

\begin{thm}
 If $L$ is a translating solution that satisfies (H), is almost-calibrated, and static, then $L$ is a plane.
\end{thm}
\begin{proof}
	In view of Lemma \ref{decomp2} and Theorem \ref{main}, it suffices to show
	\begin{prop}\label{technical} Let $L$ be a translating solutions satisfying hypothesis (H), almost-calibrated, and static. Then for  almost all $-\pi/2<a<\pi/2$ we have
				$$\{\theta=a\} \subset B_{R(a)}	$$
		for some positive $R(a)>0$.  
\end{prop}
\begin{proof}
The static condition implies the existence of a sequence of blow-downs
			$$L^i_s:=\lambda_i L_{s/\lambda_i^2}\quad\mbox{where}\quad \lim_{i\to\infty}\lambda_i=0$$
			converging weakly for every $s$ to 
			$$L^{\infty}=m_1P_1+\cdots+m_NP_N,$$
			where  $P_1,\ldots,P_N$ are Lagrangian planes with multiplicities $m_1,\ldots,m_N$ and Lagrangian angles $\bar\theta_1,\ldots,\bar\theta_N$.
			
			We claim that 
			$$P_j:=\{(u,0,v\cos \bar \theta_j ,v\sin \bar \theta_j)\,|\, (u,v)\in \R^2\}.$$			The reason is that the coordinate $y_1$ equals, up to a constant,  $-\theta$  (Proposition \ref{equacoes} (i)) and so is bounded for every $\L_t$. Thus each plane $P_j$ must have $e_1^{\bot}=0$ because, for every $R<0$ and $s\leq 0$,
		$$\lim_{i\to\infty}\sup\{\,y_1\,|\, x\in L^i_s\cap B_R(0)\}=0.$$
		 Hence the almost-calibrated condition implies that  $P_j$ is  uniquely determined by its Lagrangian angle.

 From Theorem \ref{scale} we know that the set of limiting Lagrangian angles does not depend on the sequence of  rescalings chosen and so any other convergent sequence of blow-downs  $ (\bar L^i_s)$ must  converge to $L^{\infty}$ for every $s\leq 0$.  This observation will be used later.
 
 Sard's Theorem implies that, for almost all $-\pi/2<a<\pi/2$, the set $\{\theta=a\}$ is a smooth submanifold of $L$ (possible empty).We argue that only finitely many curves contained in $\{\theta=a\}$ have finite length.  Suppose that $(C_j)_{j\in\N}$ is a sequence of distinct curves contained in  $\{\theta=a\}$ having finite length. Hypothesis (H) implies that, for some $j_0$, we can find integers $b_1,\ldots b_{j_0}$ such that
			$$b_1[C_1]+\cdots+b_{j_0}[C_{j_0}]=0,$$
			where $[C_j]$ denotes the homology class of $C_j$ in $H_1(L,\Z)$. Thus, there is a compact set $K\subseteq L$ such that
			$$\partial K=b_1C_1+\cdots+b_{j_0}C_{j_0}.$$
			From Proposition \ref{equacoes} (i) we know that $\theta$ cannot have any interior maximum or minimum and  therefore, because $\theta$ equals $a$ on ${\partial K}$, $\theta$ must be constant on $K$. Analytic continuation implies that $\theta$ is constant on $L$ and this gives us a contradiction.

			Due to Huisken's monotonicity formula and  without loss of generality, we can assume that, after passing to a subsequence,
			\begin{equation*}
			\lim_{i\to\infty}\int_{L^i_{-1}} |H|^2\exp(-|x|^2)d\mu+\lim_{i\to\infty}\int_{L^i_{1}} |H|^2\exp(-|x|^2)d\mu=0.
			\end{equation*}
			Hence, from the coarea formula
		\begin{multline*}
				\int_{-\pi/2}^{\pi/2}\lambda_i\H^1(\{\theta=a\}\cap B_{\lambda_i^{-1}}(\lambda_i^{-2}e_1))da\\
				= \int_{-\pi/2}^{\pi/2}\H^1\{x\in L^i_{- 1}\cap B_1(0)\,|\,\theta^i_{ -1}=a\}da				
				 =\int_{L^i_{ -1}\cap B_{1}(0)}|H|d\mu\\
				\leq  \left(\int_{L^i_{ -1}\cap B_1(0)}|H|^2 d\mu\right)^{1/2}\left(\H^2(L^i_{ -1}\cap B_1(0)\right)^{1/2}
			\end{multline*}
		and this implies that, for almost all $-\pi/2<a<\pi/2,$
		\begin{equation*}
		\lim_{i\to\infty }\lambda_i\H^1(\{\theta=a\}\cap B_{\lambda_i^{-1}}(\lambda_i^{-2}e_1))=0.
		\end{equation*} 
		Likewise, we also obtain
			\begin{equation*}\label{thorpe2}
		\lim_{i\to\infty }\lambda_i\H^1(\{\theta=a\}\cap B_{\lambda_i^{-1}}(-\lambda_i^{-2}e_1))=0.
		\end{equation*} 
			Choose $a$ distinct from $\bar\theta_1,\cdots,\bar\theta_N$, such that 
			\begin{multline}\label{thorpe}
			\lim_{i\to\infty }\lambda_i \H^1\left(\{\theta=a\}\cap B_{\lambda_i^{-1}}(\lambda_i^{-2}e_1
			)\right)\\
			+\lim_{i\to\infty}\lambda_i\H^1\left(\{\theta=a\}\cap B_{\lambda_i^{-1}}(-\lambda_i^{-2}e_1)\right)
			=0,
			\end{multline}
			and such that $\{\theta=a\}$ is a smooth submanifold of $L$.

			It suffices to show that there is no connected curve $C$ with infinite length contained in $\{\theta=a\}.$ If not, we could  find an unbounded sequence  $(x_i)$ in $C$ 
			$$x_i=t_ie_1+u_i,\quad\mbox{where}\quad |x_i|=\lambda_i^{-2}\quad\mbox{and}\quad \langle u_i,e_1\rangle=0.$$
			\begin{lemm}
			$$\liminf_{i\to\infty} |u_i|^2t_i^{-1}>0$$
			\end{lemm}
			\begin{proof}
				Suppose for some subsequence
				$$\lim{i\to\infty} |u_i|^2t_i^{-1}=0.$$
			Then
			$$\lim_{i\to\infty}|\lambda_i^{-2}-|t_i||=\lim_{i\to\infty}|t_i|\left|\sqrt{1+|u_i|^2t_i^{-2}}-1\right|\leq \lim_{i\to\infty}|t_i|^{-1}|u_i|^{2}=0$$
			and
			$$\lim_{i\to\infty}\lambda_i^2{|x_i-t_ie_1|^2}=\lim_{i\to\infty}\lambda_i^2 |t_i| |u_i|^2|t_i|^{-1}=\lim_{i\to\infty}\frac{|t_i|}{\sqrt{t_i^2+|u_i|^2}} |u_i|^2|t_i|^{-1}=0.$$		
			Thus, for every $i$ sufficiently large, $x_i$ belongs to either $B_{\lambda_i^{-1}/2}(\lambda_i^{-2}e_1)$ or $B_{\lambda_i^{-1}/2}(-\lambda_i^{-2}e_1)$ and so
			$$\H^1(C\cap B_{\lambda_i^{-1}}(\pm \lambda_i^{-2}e_1))\geq \pi\lambda_i^{-1}.$$
			This contradicts identity \eqref{thorpe}.
			\end{proof}
			After passing to a subsequence we can assume that
			$$ \lim_{i\to\infty} |t_i||u_i|^{-2}=s_1\geq 0.$$
			Moreover, we also assume without loss of generality that 
			$$\bar L^i_s:=|u_i|^{-1}L_{s |u_i|^{2}}, \quad\mbox{where}\quad -\infty<s<\infty$$
		converges  for all $s$, the sequence of manifolds $L_i:=L-x_i$ converges to a smooth translating solution $L_{\infty}$ with $\theta_{\infty}(0)=a$,  $t_i|u_i|^{-2}$ converges to $t_0\geq 0$, and $v_i:=u_i|u_i|^{-1}$ converges to a vector $v$ perpendicular  to $e_1$. The comments made at the beginning of this proof imply   $\bar L^i_s$ converges to $L^{\infty}$ for every $s\leq 0$.		
			
			Suppose $s_1=0$. For every $\beta$, $r>0$ and $s<0$
			\begin{multline*}
				\int_{L_{\infty}} (\theta_{\infty}-\beta)^2 \Phi_{0,r}(\cdot, 0)d\mu=\lim_{i\to\infty}\int_{L} (\theta-\beta)^2 \Phi_{x_i,r}(\cdot,0)d\mu\\
				=\lim_{i\to\infty}\int_{L_{-t_i}} (\theta-\beta)^2 \Phi_{u_i,r-t_i}(\cdot,-t_i)\\
				\leq \lim_{i\to\infty}\int_{L_{-t_i+s|u_i|^2} }(\theta-\beta)^2 \Phi_{u_i,r-t_i}(\cdot, -t_i+s|u_i|^2)\\
				=	\lim_{i\to\infty}\int_{\bar L^i_{-t_i|u_i|^{-2}+s} }(\theta-\beta)^2 \frac{\exp\left(-\frac{|x-v_i|^2}{4(r|u_i|^{-2}-s)}\right)}{4\pi(r|u_i|^{-2}-s)}d\mu\\
				=\int_{L^{\infty}_{s}}(\theta^{\infty}_{s-t_0}-\beta)^2\frac{\exp\left(\frac{|x-v|^2}{4s}\right)}{-4\pi s}d\mu\\	
					=\sum_{j=1}^Nn_j(\bar \theta_j-\beta)^2 \int_{P_j}\frac{\exp\left(\frac{|x-v|^2}{4s}\right)}{-4\pi s}d\mu.	
			\end{multline*}
			If $v$ did not belong to any $P_j$, we could make $s$ go to zero to conclude that the leftmost hand-side of the inequalities above is zero for every $\beta$, which is impossible. The fact that $v$ is perpendicular to $e_1$ implies that it can  belong at most to  one $P_{j_0}$ and thus, making $r$ go to zero and then $s$ go to zero, we see that 
			$$(a-\beta)^2\leq m_{j_0}(\bar\theta_{j_0}-\beta)^2$$
			for every $\beta$.  Therefore, $a=\bar\theta_{j_0}$ and so $s_1$ must be positive.
			
                    If $s_1>0$, then
			$$\bar L_s^i= |u_i|^{-1}L_{s|u_i|^2}=\lambda_i^{-1}|u_i|^{-1}L^i_{s|u_i|^2\lambda_i^{2}}$$
			converges weakly to
			$\sqrt{ s_1}L^{\infty}_{s/s_1}=L^{\infty}$  {for every}  $ -\infty<s<\infty$
			because
			$$\lim_{i\to\infty} |u_i|^{-2}\lambda_i^{-2}=\lim_{i\to\infty} \sqrt{t_i^{2}|u_i|^{-4}+|u_i|^{-2}}=s_1.$$
			As a result, we can argue as in the $s_1=0$ case and get a contradiction with the way $a$ was chosen.

\end{proof}
\end{proof}

The proof of Theorem A will be completed after we show
\begin{thm}
If $L$ is a translating solution that satisfies hypothesis (H) and has
$$\int_L |H|^2 d\mu \leq C_2$$ 
for some $C_2$, then $L$ is a plane.
\end{thm}
\begin{proof}

In light of Lemma \ref{decomp2} and Theorem \ref{main}, it suffices to show that for  almost all $-C_1<a<C_1$ we have
				$$\{\theta=a\} \subset B_{R(a)}	$$
		for some positive $R(a)>0$.  
	Moreover,  we can argue as	 in Proposition \ref{technical} and see that the fact that $L$ has finite first Betti number implies that it is sufficient to show that for almost all $-C_2\leq s\leq C_2$ there is no curve contained in $\{\theta=s\}$ having infinite length. 
				
				Suppose that such  a smooth curve exists and denote it by $C$. Then, for all $r$ sufficiently large, $C\cap\{|x|=r\}$ is not empty and so 
			\begin{equation}\label{lbound}
				\H^1(C\cap \{r\leq |x|\leq 2r\})\geq r.
			\end{equation}
			On the other hand, using the coarea formula, we have for all $r$
			\begin{multline*}
				\int_{-C_1}^{C_1}\H^1(\{\theta=s\}\cap \{r\leq |x|\leq 2r\})ds=  \int_{L\cap \{r\leq |x|\leq 2r\}}|H|d\mu\\
				\leq  \left(\int_{L\cap \{r\leq |x|\leq 2r\}}|H|^2 d\mu\right)^{1/2}\left(\H^2(L\cap \{r\leq |x|\leq 2r\})\right)^{1/2}\\
				\leq  \sqrt C_1 r \left(\int_{L\cap \{r\leq |x|\leq 2r\}}|H|^2 d\mu\right)^{1/2}.
			\end{multline*}
			Note that
			$$\lim_{r\to\infty}\int_{L\cap \{r\leq |x|\leq 2r\}}|H|^2 d\mu=0$$
			and thus, for almost all $-C_2\leq s\leq C_2$, we can  find a sequence $r_i$ going to infinity such that
			$$\lim_{i\to\infty} r_i^{-1}\H^1(\{\theta=s\}\cap \{r_i\leq |x|\leq 2r_i\})=0.$$
			This contradicts \eqref{lbound}.
\end{proof}

\appendix

\section{ }
Let $\tau$ be the volume form of a plane $P$ extended by parallel translation to all of $\C^2$. We denote by $\p_P$ the projection onto the plane $P$ and  assume that $P$ contains the line spanned by $e_1$.
Given any surface $\Sigma$ in $\C^2$, we denote by $\ast\tau$ the Hodge-dual of $\tau$ and if $F$ denotes an immersion of $\Sigma$ on $\C^2$, we also use $\p_P$ to represent $\p_P\circ F$. 

In what follows, $\Sigma$ will be a complete noncompact surface with smooth boundary such that $\widehat B_{r_0}(x)$ is simply-connected for all $x$ in $L$,$$\sup_{R>0} R^{-2}\H^2(\Sigma\cap B_R(0))\leq D\quad\mbox{and}\quad\inf_{x\in\Sigma,\, r\leq r_0} r^{-2}\H^2(\widehat B_{r}(x))\geq D^{-1}$$
for some constant $D$.

 \begin{lemm}\label{basic2}
 	Assume that we can find $R>0$ and $\varepsilon>0$ such that
 	$$\partial \Sigma\subset B_R(0)$$
 	and $$\inf_{\widehat B_{r_0}(x)}\ast\tau\geq \varepsilon\quad\mbox{for all} \quad x\in \Sigma\setminus B_R(0).$$
	
	Outside a compact set, $\Sigma$ can be decomposed into $k$ connected components $\Sigma_j$, $j=1,\ldots,N$ having the following property.
	\begin{itemize}
			\item[(i)] For every $x$ in $\Sigma\setminus B_R(0)$, $\widehat B_{r_0}(x)$ can be written as the graph of a function defined over $P$ with its derivatives  bounded by a constant depending only on $\varepsilon$;			
			\item[(ii)] $$\p_P\colon \Sigma_j \longrightarrow P\setminus B_{R}(0)$$ is a $n_j$-fold covering map.		
	\end{itemize}
 \end{lemm}
	\begin{proof}
		The first property is an immediate consequence of the fact that $\widehat B_{r_0}(x)$ is simply-connected. Consider the map 
		$$\p_P\colon \Sigma\setminus B_{R}(0) \longrightarrow P.$$
		The local graphical property combined with the uniform lower bounds on area densities implies that $\p_P^{-1}(B_R(0))$ is a compact subset of  $\Sigma\setminus B_{R}(0)$. 
		
		Decompose $\p_P^{-1}(P\setminus B_{R}(0))$ into connected components $(\Sigma_j)_{j\in \N}$. The local graphical property  implies the existence of an integer $n_j$ so that 
		$$\p_P\colon \Sigma_j \longrightarrow P\setminus B_{R}(0)$$
		is a $n_j$-fold covering map with
		$$\p_P(\partial \Sigma_j)\subseteq \partial (P\setminus B_R(0)).$$
		Moreover, there is a constant $C=C(\varepsilon)$ such that
		$$\lim_{R\to\infty} R^{-2}\H^2(\Sigma_j\cap B_R(0))\geq C$$
		and so there can only exist finitely many connected components $\Sigma_j$.
	\end{proof}

\section{ }

The next proposition was proven in \cite[Appendix A]{neves} with slightly different hypothesis. For that reason we will only indicate the modifications in the proof.

\begin{prop}\label{poincare}
Let $(N^i)$ and $(\alpha_i)$ be a sequence of smooth Lagrangian surfaces in $\R^4$  and  smooth functions on $N^i$
respectively, such that $(N^i)$ converges as Radon measure and as currents to a union of planes with positive integer multiplicities $N$. We
assume that, for some $R>0$, the following properties hold:
\begin{itemize}
\item[a)] There exists a constant $D_0$ such that  $$\H^2(N^i\cap B_{3R}))\leq D_0R^2$$ and
$$\cos \theta^i\geq D_0^{-1}$$for all $i\in\N$.

\item[b)]$$\lim_{i \to \infty}\int_{N^i\cap B_{3R}(0)}|\nabla \alpha_i|^2d\mu=0.$$

\item[c)] There exists a constant $D_1$ for which
$$\sup_{N^i\cap B_{3R}(0)}|\nabla \alpha_i|+R^{-1}\sup_{N^i\cap B_{3}(0)}|\alpha_i|\leq D_1$$ for all $i \in \N$.

\item[d)] For all $i \in\N$, $$N^i\cap B_{2R}(0)\quad\mbox{ is connected}$$ and
$$ \partial(N^i\cap B_{3R}(0))\subset\partial B_{3R}(0).$$
\end{itemize}
Then, there is  a real number $\alpha$ such that, after passing to a subsequence, we have for all $\phi$ with
compact support in $B_R(0)$ and all $f$ in $C(\R)$
 $$\lim_{i \to
   \infty}\int_{N^i}f(\alpha_i)\phi d\mu=f(\alpha)\mu_N(\phi),$$
where $\mu_N$ denotes the Radon measure associated to $N$.
\end{prop}

\begin{proof} It suffices to
find $\alpha \in \R$ and a sequence $(\varepsilon_j)$ converging to zero such that, for some appropriate
subsequence, we have for all $j \in \N$
$$\lim_{i \to \infty}\H^2(\{|\alpha_i-\alpha|\leq \varepsilon_j\}\cap B_{R}(0))=\H^2(N\cap B_R(0)).$$

For the rest of this proof, $K=K(D_0,D_1,k)$ will denote a generic constant depending only on the mentioned
quantities. Choose any sequence $(x_i)$ in $N^i\cap B_R(0)$. After passing to a subsequence, we have that $$
\lim_{i \to\infty}x_i=x_0\quad\mbox{and}\quad\lim_{i \to
  \infty}\alpha_i(x_i)=\alpha
$$ for some $x_0 \in B_R(0)$ and $\alpha \in \R$. Furthermore, consider also
a sequence $(\varepsilon_j)$ converging to zero such that, for all $j\in \N$,
$$
\lim_{i \to \infty}\H^{1}\bigl( \{\alpha_i=\alpha\pm\varepsilon_j\}\cap B_{3R}\bigr)=0.
$$
Such a subsequence exists because, by the coarea formula, we have
\begin{multline*}
\lim_{i \to \infty}\int_{-\infty}^{\infty}\H^{1}\bigl( \{\alpha_i=s\}\cap B_{3R}\bigr)ds
=\lim_{i \to \infty}\int_{N^i\cap B_{3R}}|\nabla \alpha_i|d\mu\\
\leq \lim_{i \to \infty}K R\left(\int_{N^i\cap B_{3R}}|\nabla \alpha_i|^2d\mu\right)^{1/2} =0.
\end{multline*}
Define 
$$
 N^{i,\alpha,j}\equiv\{|\alpha_i-\alpha|\leq \varepsilon_j\}.
 $$

Standard compactness theorems imply that, after passing to a subsequence, we have
convergence to a boundaryless integral current $\bar N^{\alpha,j}$  and to a Radon measure $ N^{\alpha,j}$,
 both having their support contained in $N$. The almost-calibrated condition implies that the support of  $N^{\alpha,j}$ coincides with the support of $\bar N^{\alpha,j}$ (the argument is the same as the one used in the beginning of the proof of Theorem \ref{scale}). The Constancy Theorem \cite[Theorem 26.27]{Leon} implies  that the support of $N^{\alpha,j}$ is a union of planes with multiplicities.

We note that  $N^i$ being almost calibrated  (see \cite[Lemma 7.1]{neves}) implies the existence of some constant $D$ such that
	$$\left(\H^2(A)\right)^{1/2}\leq D \H^{1}(\partial A),$$
where $A$ is any open subset of $N$ with rectifiable boundary. The rest of the proof follows exactly like it was done in \cite[Proposition A. 1]{neves}.
\end{proof}

\newpage
\bibliographystyle{amsbook}

\begin{thebibliography}{99}




\bibitem {huisken} G. Huisken,
Asymptotic behavior for singularities of the mean curvature flow, {\bf J. Differential Geom. 31} (1990),
285--299.
\bibitem {ilmanen1} T. Ilmanen, Elliptic Regularization and Partial Regularity for Motion by Mean Curvature. {\bf Mem. Amer. Math. Soc. 108} (1994),  1994.
\bibitem {ilmanen} T. Ilmanen, Singularities of Mean Curvature Flow of Surfaces. Preprint.
\bibitem {JYT}Y--I Lee, Private communication.
\bibitem {neves} A. Neves, Singularities of Lagrangian Mean Curvature Flow: Zero-Maslov class case. {\bf Invent. Math. 168} (2007), 449--484.
\bibitem{perelman} G. Perelman, The entropy formula for the Ricci flow and its geometric applications. preprint.

\bibitem {Leon} L. Simon,
Lectures on geometric measure theory, {\bf Proceedings of the Centre for Mathematical Analysis, Australian
National University, 3}.
\bibitem{smoczyk} K. Smoczyk, Existence of solitons for the Lagrangian mean curvature flow. 
ETH Z\"urich; 1998; preprint.
\bibitem{white1} B. White, The size of the singular set in mean curvature flow of mean-convex sets. 
{\bf J. Amer. Math. Soc. 13} (2000), 665--695. 
\bibitem{white2} B. White, The nature of singularities in mean curvature flow of mean-convex sets. 
{\bf J. Amer. Math. Soc. 16} (2003), 123--138 .

\end{thebibliography}

\vspace{20mm}

\end{document}